\documentclass[12pt]{amsart}  
\usepackage{amssymb}
\usepackage{amsmath}

\newcommand{\lqs}{\leqslant}
\newcommand{\gqs}{\geqslant}

\begin{document}

{\small Izv. Inst. Mat. i Inform. Udmurt. Univ., vyp.~3 (37), Izhevsk, 2006, 
pp. 27-28 (in Russian)}

\vskip 2.0cm

\title[Weyl almost periodic selections]{Weyl almost periodic selections of 
multivalued maps}
\author{L.I.{\,}Danilov}
\date{}

\address{Physical-Technical Institute, Ural Branch of the Russian Academy of Sciences,
Kirov Street 132, Izhevsk 426000, Russia}
\email{danilov@otf.pti.udm.ru}



\maketitle

Let $(U,\rho )$ be a complete metric space and let a point $x_0\in U$ be fixed. For
$p \gqs 1$, we denote by $M_p({\mathbb R},U)$ the set of (strongly) measurable
functions $f:{\mathbb R}\to U$ such that
$$
\sup\limits_{\xi \in {\mathbb R}}\, \int\limits_{\xi}^{\xi +1}\rho ^{\, p}\, (f(t),
x_0)\, dt < +\infty \, .
$$
On $M_p({\mathbb R},U)$ we introduce the metrics
$$
D^{(\rho )}_{p,\, l}(f,g)= \Bigl( \, \sup\limits_{\xi \in {\mathbb R}}\, \frac 1l\, 
\int\limits_{\xi}^{\xi +l}\rho ^{\, p}\, (f(t),g(t))\, dt \Bigr) ^{1/p},\ l>0\, ,
$$
and the semimetric 
$$
D^{(\rho )}_p(f,g)=\lim\limits_{l\to +\infty }\, D^{(\rho )}_{p,\, 
l}(f,g)\, ,\ f,g\in M_p({\mathbb R},U)\, .
$$
For a measurable set $T\subseteq {\mathbb R}$ (in the sense of Lebesgue), we use
the notation
$$
\varkappa _W(T)=\lim\limits_{l\to +\infty }\, \sup\limits_{\xi \in {\mathbb R}}\
\frac 1l\ {\rm {meas}}\, [\xi ,\xi +l]\cap T\, ,
$$
where $meas$ is Lebesgue measure on ${\mathbb R}$. Let $M^*_p({\mathbb R},U)$ be the
set of functions $f\in M_p({\mathbb R},U)$ such that for any $\varepsilon >0$ 
there is a number $\delta >0$ such that for all measurable sets $T\subseteq {\mathbb 
R}$ with $\varkappa _W(T)<\delta $ the inequality
$$
\lim\limits_{l\to +\infty }\, \sup\limits_{\xi \in {\mathbb R}}\
\frac 1l\ \int\limits_{[\xi ,\xi +l]\, \cap \, T}\rho ^{\, p}\, (f(t),x_0)\, 
dt < \varepsilon 
$$ 
holds.

Let $\varepsilon >0$. For a function $f\in M_p({\mathbb R},U)$, a number $\tau 
\in {\mathbb R}$ is called {\it $(\varepsilon ,D^{(\rho )}_{p,\, l})$-almost period}
if $D^{(\rho )}_{p,\, l}(f(.),f(.+\tau ))<\varepsilon $, and is called {\it 
$(\varepsilon ,D^{(\rho )}_p)$-almost period} if $D^{(\rho )}_p(f(.),f(.+\tau ))
<\varepsilon $. A set $T\subseteq {\mathbb R}$ is {\it relatively dense} if there
exists a number $a>0$ such that $[\xi ,\xi +a]\cap T\neq \emptyset $ for all 
$\xi \in {\mathbb R}$. A function $f\in M_p({\mathbb R},U)$, $p \gqs 1$, belongs 
to the space $W_p({\mathbb R},U)$ of {\it Weyl almost periodic (a.p.)} functions 
of {\it order $p$} if for any $\varepsilon >0$ there is a number $l=l(\varepsilon 
,f)>0$ such that the set of $(\varepsilon ,D^{(\rho )}_{p,\, l})$-almost periods of
$f$ is relatively dense. 

On the space $U$, we also define the metric $\rho ^{\, \prime}(x,y)=\min \{ 1,\rho
(x,y)\} $, $x,y\in U$. Let $W({\mathbb R},U)\doteq W_1({\mathbb R},(U,\rho 
^{\, \prime}))$ be the space of {\it Weyl a.p.} functions $f:{\mathbb R}\to U$
(of order 1) taking values in the metric space $(U,\rho ^{\, \prime})$. We have
$$
W_p({\mathbb R},U)\subseteq W_1({\mathbb R},U)\subseteq W({\mathbb R},U)\, ,\
p \gqs 1\, .
$$
A sequence $\{ \tau _j\} _{j\, \in \, {\mathbb N}}\subset {\mathbb R}$ is said to be 
{\it $f$-returning} for a function $f\in W({\mathbb R},U)$ if 
$D^{(\rho ^{\, \prime})}_1(f(.),f(.+\tau _j))\to 0$ as $j\to +\infty $. For a function
$f\in W({\mathbb R},U)$, by ${\rm {Mod}}\, f$ we denote the module (additive group) 
of all numbers $\lambda \in {\mathbb R}$ such that $e^{\, {\rm i}\lambda \tau
_j}\to 1$ as $j\to +\infty $ (where ${\rm i}^2=-1$) for all $f$-returning 
sequences $\{ \tau _j\} $. If $U=({\mathcal H},\| .\| )$ is a Banach space and 
$f\in W_1({\mathbb R},{\mathcal H})$, then ${\rm {Mod}}\, f$ coincides with the
frequency module of the function $f$.

Let $({\rm {cl}}_{\, b}\, U, {\rm {dist}}_{\rho})$ be the metric space of non-empty
closed bounded sets of the space $(U,\rho )$ with the Hausdorff metric ${\rm {dist}}
_{\rho}\, $; $({\rm {cl}}\, U, {\rm {dist}}_{\rho ^{\, \prime}})$ is the metric 
space of non-empty closed sets of the space  $(U,\rho )$ with the Hausdorff metric
${\rm {dist}}_{\rho ^{\, \prime}}\, $ corresponding to the metric $\rho ^{\, \prime}$. 
The spaces $W_p({\mathbb R},{\rm {cl}}_{\, b}\, U)$, $p \gqs 1$, and $W({\mathbb R},
{\rm {cl}}_{\, b}\, U)$ of {\it Weyl a.p. multivalued maps} $F:{\mathbb R}\to {\rm 
{cl}}_{\, b}\, U$ are defined as the spaces of Weyl a.p. functions taking values in 
the metric space $({\rm {cl}}_{\, b}\, U, {\rm {dist}}_{\rho})$. Let $W({\mathbb R},
{\rm {cl}}\, U)\doteq W_1({\mathbb R},({\rm {cl}}\, U,{\rm {dist}}_{\rho ^{\, 
\prime}}))$; $W({\mathbb R},{\rm {cl}}_{\, b}\, U)\subseteq W({\mathbb R},{\rm {cl}}\, 
U)$.
\vskip 0.2cm

{\bf Theorem 1} (see \cite{1}). {\it 
Let $(U,\rho )$ be a complete metric space, let $g\in W({\mathbb R},U)$, and let
$F\in W({\mathbb R},{\rm {cl}}\, U)$. Then for any increasing function $\eta :
[0,+\infty )\to [0,+\infty )$, for which $\eta (0)=0$ and $\eta (t)>0$ for all $t>0$,
there is a function $f\in W({\mathbb R},U)$ such that ${\rm {Mod}}\, f\subseteq
{\rm {Mod}}\, g+{\rm {Mod}}\, F$, $f(t)\in F(t)$ a.e., and 
$$
\rho (f(t),g(t))\lqs \rho (g(t),F(t))+\eta (\rho (g(t),F(t)))\ a.e.
$$ 
If $F\in W_p({\mathbb R},{\rm {cl}}_{\, b}\, U)$ for some $p\gqs 1$, then the 
function $f$ belongs to the space $W_p({\mathbb R},U)$.}
\vskip 0.2cm

For a complete separable metric space $(U,\rho )$ we denote by $({\mathcal M}(U),d)$ 
the (complete separable) metric space of Borel probability measures $\mu [.]$
defined on the $\sigma$-algebra of Borel subsets of the metric space $(U,\rho )$, 
with the L${\acute {\mathrm e}}$vy -- Prokhorov metric $d$; $W_1({\mathbb R},{\mathcal 
M}(U))$ is the space of {\it Weyl a.p. measure-valued functions} ${\mathbb R}\ni t\to 
\mu [.;t]\in {\mathcal M}(U)$ (of order 1) taking values in the metric space 
$({\mathcal M}(U),d)$. A  measure-valued function ${\mathbb R}\ni t\to \mu [.;t]\in 
{\mathcal M}(U)$ belongs to the space $W_1({\mathbb R},{\mathcal M}(U))$ if and only if 
$$
\int\limits_U{\mathcal F}(x)\, \mu [dx;.]\in W_1({\mathbb R},{\mathbb R})
$$ 
for all continuous bounded functions ${\mathcal F}\in C_b(U,{\mathbb R})$ (see 
\cite{2}). Moreover,
$$
{\rm {Mod}}\, \mu [.;.]=\sum\limits_{{\mathcal F}\, \in \, C_b(U,{\mathbb R})}{\rm 
{Mod}}\, \int\limits_U{\mathcal F}(x)\, \mu [dx;.]\, .
$$
For $\mu \in {\mathcal M}(U)$, $x\in U$, and $\delta \in (0,1)$ we shall use the
notation 
$$
r_{\delta}(x,\mu )=\inf \, \{ r>0:\mu [U_r(x)]>\delta \} \, ,
$$
where $U_r(x)=\{ y\in U:\rho (x,y)<r\} $.
\vskip 0.2cm

{\bf Theorem 2}. {\it
Let $(U,\rho )$ be a complete separable metric space, let $g\in W({\mathbb R},U)$, 
and let $\mu [.;.]\in W_1({\mathbb R},{\mathcal M}(U))$. Then for any $\delta \in 
(0,1)$ there exists a function $f_{\delta}\in W({\mathbb R},U)$ such that ${\rm 
{Mod}}\, f_{\delta}(.)\subseteq {\rm {Mod}}\, g(.)+{\rm {Mod}}\, \mu [.;.]$, 
$f_{\delta}(t)\in {\rm {supp}}\, \mu [.;t]$ a.e., and 
$$
\rho (f_{\delta}(t),g(t))< r_{\delta}(g(t),\mu [.;t])+\delta \ a.e.
$$ 
}
\vskip 0.2cm

Let ${\mathcal S}_{\, \rm {rd}}$ be the collection of relatively dense subsets 
$T\subseteq {\mathbb R}$. We shall denote by ${\mathcal P}^{(\rho )}_p
(\varepsilon ;f)$ the set of $(\varepsilon ,D^{(\rho )}_p)$-almost periods of a 
function $f\in M_p({\mathbb R},U)$. Let $\widetilde W_p({\mathbb R},U)$ be the set of
functions $f\in M_p({\mathbb R},U)$ satisfying the following two conditions:

1) ${\mathcal P}^{(\rho )}_p(\varepsilon ;f)\in {\mathcal S}_{\, \rm {rd}}$ for
all $\varepsilon >0$;

2) for any $\varepsilon ,\delta >0$ there is a compact set $K_{\varepsilon ,\delta }
\subseteq U$ such that $\varkappa _W(\{ t\in {\mathbb R}: \rho (f(t),K_{\varepsilon 
,\delta })\gqs \varepsilon \} ) <\delta $. (In the case $U=({\mathbb R}^n,|.|)$,
the condition 2 is satisfied for all functions $f\in M_1({\mathbb R},{\mathbb R}^n)$.) 

Let $\widetilde W({\mathbb R},U)\doteq \widetilde W_1({\mathbb R},(U,\rho 
^{\, \prime}))$. We have 
$$
W_p({\mathbb R},U)=W({\mathbb R},U)\cap M^*_p({\mathbb R},U)\subseteq \widetilde 
W({\mathbb R},U)\cap M^*_p({\mathbb R},U)\subseteq \widetilde W_p({\mathbb R},U)\, .
$$
If $f\in \widetilde W({\mathbb R},U)\cap M^*_p({\mathbb R},U)$, then for any 
$\varepsilon ^{\, \prime}>0$ there exists a number $\varepsilon >0$ such that
${\mathcal P}^{(\rho )}_p(\varepsilon ^{\, \prime};f)\supseteq {\mathcal P}^{(\rho 
^{\, \prime})}_1(\varepsilon ;f)$.
\vskip 0.2cm

{\bf Theorem 3}. {\it
Let $(U,\rho )$ be a complete metric space, let $g\in \widetilde W({\mathbb R},U)$, 
and let $F\in \widetilde W_1({\mathbb R},({\rm {cl}}\, U,{\rm {dist}}
_{\rho ^{\, \prime}}))$.  Suppose that 
$$
{\mathcal T}(\varepsilon )\doteq {\mathcal 
P}^{(\rho ^{\, \prime})}_1(\varepsilon ;g)\cap {\mathcal P}^{({\rm {dist}}_{\rho 
^{\, \prime}})}_1(\varepsilon ;F)\in {\mathcal S}_{\, \rm {rd}}
$$ 
for all $\varepsilon >0$. Then for any $\delta >0$ and any function $h\in W_1
({\mathbb R},{\mathbb R})$ with a dense {\rm (}countable{\rm )} 
frequency module ${\rm {Mod}}\, h$ {\rm (}for this function ${\mathcal T}
(\varepsilon )\cap {\mathcal P}^{(\rho _{\, \mathbb R})}_1(\varepsilon ;h)\in 
{\mathcal S}_{\, \rm {rd}}$ for all $\varepsilon >0${\rm )} there exists a function 
$f_{\delta}\in \widetilde W({\mathbb R},U)$ such that $f_{\delta}(t)\in F(t)$ a.e., 
$$
\rho (f_{\delta}(t),g(t))< \rho (g(t),F(t))+\delta \ a.e.\, ,
$$ 
and for any $\varepsilon ^{\, \prime}>0$ there is a number $\varepsilon >0$
such that 
$$
{\mathcal P}^{(\rho ^{\, \prime})}_1(\varepsilon ^{\, \prime};f_{\delta})\supseteq 
{\mathcal T}(\varepsilon )\cap {\mathcal P}^{(\rho _{\, \mathbb R})}_1(\varepsilon 
;h)\, .
$$
Furthermore, if $F\in M^*_p({\mathbb R},{\rm {cl}}_{\, b}\, U)$ for some $p \gqs 1$, 
then $f_{\delta}\in \widetilde W({\mathbb R},U)\cap M^*_p({\mathbb R},U)
\subseteq \widetilde W_p({\mathbb R},U)$.}

\vskip 2.0cm

\end{document}